\documentclass[11pt]{article}
\usepackage{amsmath,amssymb,amsthm}
\usepackage[margin=1in]{geometry}
\usepackage[hidelinks]{hyperref}
\newtheorem{theorem}{Theorem}
\newtheorem{lemma}[theorem]{Lemma}
\newtheorem{corollary}[theorem]{Corollary}
\newcommand{\ZZ}{\mathbb Z}

\newcommand{\SL}{\mathrm{SL}}
\title{A regular free action on the integral plane and a problem on left braces}
\author{L\"u Gong\\[3pt]
  \small School of Mathematics and Statistics, Nantong University\\
  \small Jiangsu 226019, China\\
  \small\href{mailto:lieningzai1917@126.com}{\nolinkurl{lieningzai1917@126.com}}}\date{}
\begin{document}
\maketitle
\begin{abstract}
We construct a free subgroup of countably infinite rank in
$\ZZ^2\rtimes\SL_2(\ZZ)$ that acts regularly on $\ZZ^2$.
At each stage we adjoin an affine transformation taking the origin to
the least point of a prescribed orbit. Disjoint domains for reduced words
ensure that the stabilizer of the origin remains trivial.
The orbit map gives a left brace with additive group $\ZZ^2$ and
multiplicative group free of countably infinite rank. This group contains
a nonabelian free subgroup and is therefore nonsolvable, answering both
parts of Problem 20.126 of the Kourovka Notebook.
\end{abstract}



\section{Introduction}
Bardakov, Neshchadim and Yadav ask in Problem 20.126 of the
Kourovka Notebook \cite{KN} whether a left brace with finitely
generated abelian additive group can have a nonsolvable multiplicative
group, and whether its multiplicative group can contain a nonabelian free
subgroup. A left brace has an abelian group law $+$ and a group law
$\circ$ satisfying $x\circ(y+z)+x=(x\circ y)+(x\circ z)$.
A single left brace answers both questions.

\begin{theorem}\label{main}
The group $\ZZ^2\rtimes\SL_2(\ZZ)$ contains a free subgroup $\Gamma$ of
countably infinite rank acting regularly on $\ZZ^2$. Consequently there
is a left brace with additive group $\ZZ^2$ and multiplicative group
isomorphic to $\Gamma$.
\end{theorem}

Sat\^o \cite{Sato} constructed a free group of rank two acting freely
on the integral plane by determinant-one affine transformations. In his
construction the generators fix a common rational point $a/q$, with $a\in\ZZ^2$ and $q$ a positive integer.
Such an action has infinitely many integral orbits: the greatest common
divisor of the coordinates of $qx-a$ is constant along each orbit, since
$qg(x)-a$ is obtained from $qx-a$ by a unimodular matrix. Starting at
$x$ with $qx-a\ne0$, the integral points $nx+(1-n)a/q$ for positive
$n\equiv1\pmod q$ give infinitely many distinct values.
For a brace we need a regular action: every lattice point must be reached
from zero by exactly one group element.

We enforce uniqueness by disjoint domains for reduced words. To obtain
existence, we enlarge the group at each step by a generator taking zero
to the least point of a prescribed missing orbit. The finite-ball
assertion in Lemma~\ref{local} allows us to prescribe this image.

\section{Domains for affine transformations}
Write $X=\ZZ^2$. Order $X$ lexicographically by the triples
$(\|x\|^2,x_1,x_2)$, using the Euclidean norm, and write $\preceq$ for
this order. Every nonempty subset has a least element, since bounded
balls contain finitely many lattice points. The least point of $X$ is
zero. For a permutation $s$ of $X$ put
\[
 D(s)=\{x:s^{-1}x\preceq x\}.
\]
If $x\notin D(s^{-1})$, then $x\prec sx$ and hence $sx\in D(s)$.
Suppose that for a family of generators the sets $D(s)$ for all signed
generator symbols are pairwise disjoint and omit zero. A nonempty
reduced word sends zero into the set for its leftmost letter: read the
word from right to left and apply the preceding implication. It follows
that the generators freely generate their group and that the stabilizer
of zero is trivial. The same argument sends any point outside all the
sets into one of them under every nonempty reduced word.

For a nonzero real vector $x$ let $[x]$ denote its projective direction.
Choose integral hyperbolic matrices $L_j$, $j\ge1$, and open intervals
$I_j^+,I_j^-$ in the real projective line around their expanding and
contracting directions respectively, with pairwise disjoint closures
avoiding the vertical direction. Explicitly, take
\[
 A=\begin{pmatrix}2&1\\1&1\end{pmatrix},\quad
 T=\begin{pmatrix}1&1\\0&1\end{pmatrix},\quad L_j=T^{10j}AT^{-10j}.
\]
In the coordinate $x_1/x_2$ the two directions of $L_j$ are
$10j+(1\pm\sqrt5)/2$. Intervals of radius $1/4$ suffice. Both
eigenlines contain no nonzero rational vector.

\begin{lemma}\label{local}
Fix $j$, write $L=L_j$, and let $0\ne y\in X$. Put $g_n x=y+L^n x$.
There are $R,N$ such that, for $n\ge N$ and $\|x\|>R$,
\[
 x\in D(g_n)\ \Longrightarrow\ [x]\in I_j^+,\qquad
 x\in D(g_n^{-1})\ \Longrightarrow\ [x]\in I_j^-.
\]
Moreover, for each $B=\{x\in X:\|x\|\le r\}$ there is an integer
$N_B$ such that, for all $n\ge N_B$,
\[
 D(g_n)\cap B=\{y\}\cap B,\qquad D(g_n^{-1})\cap B=\varnothing.
\]
\end{lemma}
\begin{proof}
Let $P_+,P_-$ denote the projections onto the two eigenlines along
one another, and let $\lambda>1$ be the expanding eigenvalue. If
$x\in D(g_n)$ then $\|L^{-n}(x-y)\|\le\|x\|$. Projecting onto the
contracting eigenline gives, for $x\ne0$,
\[
 \frac{\|P_-x\|}{\|x\|}
 \le \|P_-\|\lambda^{-n}+\frac{\|P_-y\|}{\|x\|}.
\]
The continuous function $u\mapsto\|P_-u\|$ has a positive minimum
on the compact set of unit vectors with $[u]\notin I_j^+$, since
$\ker P_-$ is the expanding eigenline. The bound
is smaller than this minimum once both $n$ and $\|x\|$ exceed fixed
thresholds. Likewise $x\in D(g_n^{-1})$ gives
$\|y+L^nx\|\le\|x\|$ and hence
\[
 \frac{\|P_+x\|}{\|x\|}
 \le\|P_+\|\lambda^{-n}\left(1+\frac{\|y\|}{\|x\|}\right).
\]
Use the positive minimum off $I_j^-$ and take common thresholds.
These bounds depend on $j,y$ only.

For each fixed $x\ne y$, the nonzero integer vector $x-y$ has nonzero
contracting component by irrationality, so
$\|g_n^{-1}x\|=\|L^{-n}(x-y)\|\to\infty$.
For each fixed $x\ne0$, $\|g_nx\|\to\infty$ by the same argument
using the expanding component. Finally $g_n0=y\succ0$, so zero never
belongs to $D(g_n^{-1})$. Finiteness gives simultaneous eventual
exclusion throughout the ball. The exceptional point $y$ belongs to
$D(g_n)$ for every $n$, since $g_n^{-1}y=0\prec y$.
\end{proof}

\section{Joining the orbits}
\begin{proof}[Proof of the affine assertion of Theorem~\ref{main}]
Enumerate $X$ and start with the trivial group. At each finite stage
maintain pairwise disjoint signed sets omitting zero, and require that
every point of each set outside a finite ball has direction in its
assigned interval.
For the next enumerated point $t$, do nothing if it already belongs to
the orbit of zero. Otherwise let $y$ be the least element in the current
orbit of $t$. In particular $y\ne0$. For every old signed letter $s$,
minimality gives $y\preceq s^{-1}y$. Equality is impossible: a fixed
point of $s$ belongs to both $D(s)$ and $D(s^{-1})$. Consequently $y$
is outside all old sets.

Choose the next unused index $j$. Take $R,N$ from Lemma~\ref{local}, and choose a ball $B$ about
zero of radius at least $R$ containing $y$ and all exceptions to the
old direction conditions. Outside this ball the two new sets and the old sets have
disjoint assigned intervals, so there are no new intersections. With the ball now fixed, choose $n\ge\max\{N,N_B\}$ using
the two assertions of Lemma~\ref{local}. Within the ball $D(g_n)$
contains exactly $y$ and $D(g_n^{-1})$ is empty.
The old sets omit $y$. Thus all requirements persist after adjoining
$g_n$. Since $g_n0=y$ and $t$ lies in the old orbit of $y$, the enlarged
orbit of zero contains $t$.

Every point eventually belongs to the orbit of zero. Every word uses
only finitely many letters, so the disjoint-set argument proves freeness
and triviality of the stabilizer in the union group. There must be
infinitely many adjunctions: at a finite stage all sufficiently distant
vertical lattice points lie outside the finitely many sets, whereas every
nonidentity reduced word sends zero into one of them. Thus no finite
stage is transitive. The action of the union is regular.
\end{proof}

\section{The brace and its orbit cocycle}
\begin{proof}[Completion of the proof of Theorem~\ref{main}]
Let $\Gamma$ be the group in Theorem~\ref{main}. For $x\in X$ write
$\gamma_x$ for the unique element with $\gamma_x0=x$, and $B_x$ for
its linear part. Set $x\circ y=\gamma_x y=x+B_xy$.
Regularity gives $\gamma_{x\circ y}=\gamma_x\gamma_y$, so $\circ$
transports the group law of $\Gamma$ to $X$. Additivity of $B_x$ gives
$x\circ(y+z)+x=(x\circ y)+(x\circ z)$. The multiplicative group
contains a free group of rank two and is therefore nonsolvable.
\end{proof}

\begin{corollary}
There is a free subgroup $H\le\SL_2(\ZZ)$ of countably infinite rank
whose natural action admits a bijective cocycle $d:H\to\ZZ^2$.
\end{corollary}
\begin{proof}
The kernel of the linear projection on $\Gamma$ is a normal subgroup
of translations and is therefore abelian. A nonabelian free group has
no nontrivial abelian normal subgroup: in the Cayley tree for a free basis, commuting nonidentity elements
have the same axis. Indeed, each preserves the axis of the other and
acts on it by a translation. Thus all nonidentity elements of an
abelian normal subgroup would have a common axis. Conjugation carries
axes to axes, so the whole group would preserve that line. This is impossible, since the action on vertices is transitive
and the tree is not a line. Hence the projection identifies $\Gamma$
with a subgroup $H$. Assign to each $h\in H$ the translation part of
its unique lift. Composition gives $d(hk)=d(h)+h\,d(k)$, and regularity
makes $d$ bijective.
\end{proof}


\begin{thebibliography}{1}
\bibitem{KN} E.~I. Khukhro and V.~D. Mazurov (eds.),
\emph{The Kourovka Notebook}, Problem 20.126, by V.~G. Bardakov,
M.~V. Neshchadim and M.~K. Yadav.
\bibitem{Sato} K. Sat\^o, A free group acting on $\ZZ^2$ without fixed points,
\emph{L'Enseignement Math\'ematique} (2) \textbf{45} (1999), 189--194.
\end{thebibliography}
\end{document}